\input ./style/arxiv-general.cfg
\documentclass[MSNbibl,citesort,seceqn,number,dvips]{arxbj}
\makeatletter
   \@ifpackageloaded{graphicx}{}{\usepackage{graphicx}}
\makeatother
\usepackage{mathbh}
\usepackage{upgreek}

\volume{22}
\issue{2}
\pubyear{2016}
\firstpage{1113}
\lastpage{1130}
\doi{10.3150/14-BEJ688}
\docsubty{FLA}

\makeatletter
\newcommand{\rrvert}{\vert}
\newcommand{\rrVert}{\Vert}
\newcommand{\llvert}{\vert}
\newcommand{\llVert}{\Vert}
\renewcommand{\mid}{|}
\newremark{remark}[proposition]{Remark}
\newremark{example}{Example}[section]
\newtheorem{proposition}{Proposition}[section]
\newtheorem{lemma}[proposition]{Lemma}
\newtheorem{theorem}[proposition]{Theorem}
\newtheorem{corollary}[proposition]{Corollary}
\def\la{\lambda}
\def\I{\mathbh{1}}
\def\R{{\mathbb R}}
\def\E{{\mathbb E}}
\def\P{{\mathbb P}}
\makeatother

\begin{document}
\begin{frontmatter}

\title{Excursion probability of Gaussian random fields on sphere}
\runtitle{Excursion probability of Gaussian random fields on sphere}

\begin{aug}
\author[A]{\inits{D.}\fnms{Dan}~\snm{Cheng}\thanksref{A}\ead[label=e1]{cheng@stt.msu.edu}}
\and
\author[B]{\inits{Y.}\fnms{Yimin}~\snm{Xiao}\corref{}\thanksref{B}\ead[label=e2]{xiao@stt.msu.edu}\ead[label=u2,url]{http://www.stt.msu.edu/\textasciitilde xiaoyimi}}
\address[A]{Department of Statistics, North Carolina State University,
2311 Stinson Drive, Campus Box 8203, Raleigh, NC 27695, USA. \printead{e1}}
\address[B]{Department of Statistics and Probability, Michigan State
University, 619 Red Cedar Road, C-413 Wells Hall, East Lansing, MI
48824, USA. \printead{e2};\\ \printead{u2}}
\end{aug}

%
\received{\smonth{1} \syear{2014}}
%
\revised{\smonth{6} \syear{2014}}

%
\begin{abstract}
Let $X= \{X(x)\dvt x\in\mathbb{S}^N\}$ be a real-valued, centered
Gaussian random
field indexed on the $N$-dimensional unit sphere $\mathbb{S}^N$.
Approximations to the
excursion probability $\P\{\sup_{x\in\mathbb{S}^N} X(x) \ge u
\}$, as $u\to\infty$,
are obtained for two cases: (i) $X$ is locally isotropic and its sample
functions are non-smooth
and; (ii) $X$ is isotropic and its sample functions are twice
differentiable. For case (i), the excursion probability
can be studied by applying the results in
Piterbarg (\textit{Asymptotic Methods in the Theory of {G}aussian
Processes and Fields} (1996) Amer. Math. Soc.),
Mikhaleva and Piterbarg (\textit{Theory Probab. Appl.} \textbf{41}
(1997) 367--379)
and
Chan and Lai (\textit{Ann. Probab.} \textbf{34} (2006) 80--121).
It is shown that the asymptotics of $\P\{\sup_{x\in\mathbb
{S}^N} X(x) \ge u \}$ is similar
to Pickands' approximation on the Euclidean
space which involves Pickands' constant. For case~(ii), we apply the expected
Euler characteristic method to obtain a more precise approximation such
that the error
is super-exponentially small.
\end{abstract}

%
\end{frontmatter}

\section{Introduction}\label{sec1}

Even though the characterizations of isotropic covariance functions and
variograms on spheres
were given long time ago by Schoenberg \cite{Schoenberg42} and Gangolli \cite{Gangolli},
respectively, and random
fields on the sphere were studied by Obukhov \cite{Obukhov47}, Yaglom \cite{Yaglom61} and
Jones \cite{Jones63}, it is
the applications in atmospherical sciences, geophysics, solar physics,
medical imaging and
environmental sciences (see, e.g., Genovese \textit{et~al.} \cite{GenMAWasserman04},
Oh and Li \cite{OhLi04}, Stein \cite{Stein07},
Cabella and Marinucci \cite{CabellaM09}, Tebaldi and
Sans{\'o} \cite{TebaldiS09}, Hansen \textit{et~al.} \cite{Hansen11}) that have
stimulated the recent rapid development in statistics of random fields
on the sphere. Various
new random field models have been constructed and new probabilistic and
statistical methods
have been developed. For example, Jun and Stein  \cite{JunStein07,JunStein08},
Huang, Zhang and Robeson \cite{HuangZR11},
Jun \cite{Jun11},
Hitczenko and Stein \cite{HStein12}, Ma \cite{Ma12},
Du, Ma and Li \cite{DuMa13} and
Gneiting \cite{Gneiting12} have constructed
several classes of real or vector-valued random fields on spheres;
Istas \cite{Istas05,Istas06} has
constructed spherical fractional Brownian motion (SFBM), which has
fractal sample functions, and
studied its Karhunen--Lo\`eve expansion and other properties.
Lang and Schwab \cite{LS13}
characterized sample H\"older continuity and sample differentiability
of isotropic
Gaussian random fields on the two-dimensional sphere $\mathbb S^2$ in
terms of their angular
power spectra. We refer to the recent book by Marinucci and Peccati \cite
{MarinucciPecati07} for a systematic
account on theory and statistical inferences of random fields on the
sphere $\mathbb{S}^N$,
with a view towards applications to cosmology. 

In this paper, we consider a real-valued, centered (locally) isotropic
Gaussian random field $X= \{X(x)\dvt x\in\mathbb{S}^N\}$, indexed on the
$N$-dimensional unit sphere $\mathbb{S}^N$, and investigate the
asymptotic properties of the excursion probability $\P\{\sup_{x
\in\mathbb S^{N}} X(x) \ge u \}$ as $u \to\infty$. Such
excursion probabilities are important in probability theory, statistics
and their applications. In particular, we mention that the above
excursion probability has appeared in Sun \cite{Sun91}, Park and Sun \cite{ParkSun98}
for determining the $P$-value in studying exploratory projection
pursuit and, as illustrated by Sun \cite{Sun01}, is useful for constructing
simultaneous confidence region for a function $f\dvtx \mathbb{S}^N \to\R
$. In his studies of projection-based depth functions, Zuo \cite{Zuo03} has
shown that Gaussian random fields on sphere appear as scaling limit of
sample projection median (see Theorems 3.2 and 3.3 in Zuo \cite{Zuo03}) and
the excursion probability of the limiting Gaussian field is useful for
constructing confidence regions for the true projection median (see
Remark~3.2 in Zuo \cite{Zuo03}). For further information on extreme value
theory of Gaussian random fields on Euclidean spaces or manifolds and
statistical applications, we refer to Adler and Taylor \cite{AT07},
Adler, Taylor and Worsley \cite{ATW12} and Marinucci and Peccati \cite{MarinucciPecati07}.

For studying the excursion probability of $X= \{X(x)\dvt x\in\mathbb
{S}^N\}$, we will distinguish two cases:
(i) the sample function of $X$, denoted as $ X(\cdot)$, is non-smooth
and, (ii) $ X(\cdot) \in C^2$ a.s., and to apply very different methods.
In the non-smooth case, the asymptotics of the excursion probability
$\P\{\sup_{x \in\mathbb S^{N}} X(x) \ge u \}$
as $u \to\infty$ can be studied by applying the results in Piterbarg \cite{Piterbarg96},
Mikhaleva and Piterbarg \cite{MP96} or Chan and Lai \cite{ChanL06},
which are extensions of the seminal result of Pickands \cite{Pickands69} under
various local stationarity conditions. We will make use of Theorem 2.1
in \cite{ChanL06} to prove Theorem \ref{ThmPickandsapproximationsphere} in Section~\ref{sec2}, 
and the method can also be applied to other Gaussian fields on sphere
with more complicated local covariance structures, see Section~\ref
{secSSFBM} for the example of standardized spherical fractional
Brownian motion. For the smooth case, we consider isotropic Gaussian
fields on sphere. Thanks to the special representation of covariance
function (Theorem \ref{ThmSchoenberg}), we are able to apply the
general theory of Adler and Taylor \cite{AT07} to compute the
Lipschitz--Killing curvatures induced by the field and hence derive the
approximation to the excursion probability, see Theorem \ref
{ThmEulerapproximation} and Corollary \ref{CorEulerapproximation}
below. Such
an approximation is more precise than that in Theorem \ref
{ThmPickandsapproximationsphere} for the non-smooth case and the
error is super-exponentially small.

We should mention that Mikhaleva and Piterbarg \cite{MP96} have established
asymptotic results for the excursion probability of Gaussian fields on
a finite-dimensional smooth manifold in $\mathbb R^{N+1}$. Their
theorems can be applied to obtain results similar to Theorem 2.4 below
for a Gaussian field $X$ on the sphere $\mathbb S^N$, provided $X$ is
the restriction on $\mathbb S^N$ of a Gaussian field defined on
$\mathbb R^{N+1}$. This approach is very useful, but may not be able to
deal with all locally isotropic Gaussian random fields on $\mathbb
S^N$. For instance, Huang, Zhang and Robeson \cite{HuangZR11} have recently shown
that the
restriction of some commonly used stationary isotropic covariance
functions on $\R^{N+1}$ may not be a valid covariance functions on the
sphere (when the Euclidean metric is replaced by the spherical metric).
Similarly, another method proposed\vspace*{1pt} by Ma (\cite{Ma12}, Theorem 4) to obtain
valid covariance functions on $\mathbb S^N$ from those on $\mathbb
R^{N+1}$ is only able to produce a proper subset of all covariance
functions on $\mathbb S^N$ (cf. Ma \cite{Ma12}, page~775). 
Their works have motivated us to deal with Gaussian fields on sphere
directly to establish asymptotic results for $\P\{\sup_{x \in
\mathbb S^{N}} X(x) \ge u \}$. 

Motivated by Mikhaleva and Piterbarg \cite{MP96}, as well as pointed out by
an anonymous referee, it would be interesting to study the excursion
probability for Gaussian fields over Riemannian manifolds (beyond
sphere), whose covariance functions satisfy (\ref{EqPickandscondition})
with $d(x,y)$ being the geodesic distance of $x$ and $y$.
This is beyond the scope of the present paper, but we believe that a
Pickands-type approximation similar to Theorem \ref
{ThmPickandsapproximationsphere} still holds. As pointed out by an anonymous
referee, the problem in the smooth case may be more challenging because
there is no analog of the Gegenbauer polynomials to characterize the
covariance functions of Gaussian fields over general Riemannian manifolds.

We end the \hyperref[sec1]{Introduction} with some notation. Let $\llVert\cdot\rrVert
$ and
${\langle}\cdot, \cdot\rangle$ denote, respectively, the Euclidean
norm and the inner product in $\R^{N+1}$ (or in $\R^N$, which will be
clear from the context). Denote by $d(\cdot,\cdot)$ the spherical
distance on $\mathbb{S}^N$,
that is, $d(x,y)= \arccos{\langle}x, y \rangle$, $\forall x, y\in
\mathbb{S}^N$. For two functions $f(t)$ and $g(t)$, we say $f(t) \sim
g(t)$ as $t\to t_0\in[-\infty, +\infty]$ if $\lim_{t\to t_0}
f(t)/g(t) =1$.

\section{Non-smooth Gaussian fields on sphere}\label{sec2}

We start with case (i) where the sample functions of $X= \{X(x)\dvt
x\in
\mathbb{S}^N\}$ may be non-smooth.
This case is easier and we show that the asymptotics of the excursion
probability $\P\{\sup_{x \in\mathbb S^{N}} X(x) \ge u \}$,
as $u \to\infty$, can be derived from the results in Piterbarg \cite
{Piterbarg96}, Mikhaleva and Piterbarg \cite{MP96} and Chan and Lai
\cite{ChanL06}. 

\subsection{Locally isotropic Gaussian fields on sphere}
Let $X= \{X(x)\dvt x\in\mathbb{S}^N\}$ be a centered Gaussian field with
covariance function $C$
satisfying
%
%
\begin{equation}
\label{EqPickandscondition} C(x,y)= 1-c d^\alpha(x,y) \bigl(1+\mathrm{o}(1)\bigr
)\qquad\mbox{as } d (x,y)\to0,
\end{equation}
for some constants $c>0$ and $\alpha\in(0,2]$. When $X(\cdot)$ is
smooth, we have $\alpha= 2$.

Covariance functions satisfying (\ref{EqPickandscondition}) behave
isotropically in a local sense,
hence the corresponding random fields fall under the general category
of locally isotropic random
fields. Similarly to Gaussian fields defined on the Euclidean space
(cf. Adler \cite{Adler10}), one
can show that, when $\alpha\in(0, 2)$, the sample function of $X$ is
not differentiable and the
fractal dimensions of its trajectories are determined by $\alpha$. See
Andreev and Lang \cite{AL13},
Hansen \textit{et~al.} \cite{Hansen11} and Lang and Schwab \cite{LS13} for related
regularity results.

There are many examples of covariances of isotropic Gaussian fields on
$\mathbb{S}^N$ that satisfy
(\ref{EqPickandscondition}). A~well-known example is
$C(x,y)=\mathrm{e}^{-c d^\alpha(x,y)}$, where $c>0$
and $\alpha\in(0,1]$ (cf. e.g., Huang, Zhang and Robeson \cite{HuangZR11}, page~725). In
their studies on germ-grain (or random ball)
models on the sphere $\mathbb S^N$,
Estrade and Istas (\cite{EsIstas10}, Remark 2.5 and Lemma 3.1)
discovered an
isotropic Gaussian field
$W^{\beta}$ on $\mathbb S^N$ with $0 < \beta< 1/2$, whose covariance
function satisfies
(\ref{EqPickandscondition}) for $\alpha= 2\beta\in(0,1]$. (From
here one can show that, even though
$W^\beta$ and the spherical fractional Brownian motion $B_\beta(x)$
introduced by Istas \cite{Istas05} are
different, they share some local properties (e.g., they have the same
H\"older continuity and fractal
dimensions). In Remark \ref{Re;25} below, we will compare the
excursion probabilities of $W^\beta$ and
the standardized SFBM.) Moreover, as in Yadrenko \cite{Y83} and Ma \cite{Ma12},
one can apply the identity
\[
\llVert x-y\rrVert= 2 \sin\biggl(\frac{d(x, y)} 2 \biggr)\qquad\forall x, y
\in\mathbb S^N
\]
to construct covariance functions that satisfy (\ref
{EqPickandscondition}) from isotropic covariance
functions $K(\cdot)$ on $\R^N$ which satisfy $K(x) = 1 - c_1 \llVert
x\rrVert
^\alpha(1+\mathrm{o}(1)) $ as $\llVert x\rrVert \to0$.
In particular, the following covariance function $C$ given by
Soubeyrand, Enjalbert and Sache  \cite{SES08}
%
%
\begin{equation}
\label{Defsin} C(x, y) = 1 - \biggl(\sin\frac{d(x, y)}{c^{1/\alpha}}
\biggr)^\alpha
\I_{\{d(x, y) \le\uppi c^{1/\alpha}\}},
\end{equation}
where $c > 0$ and $\alpha\in(0,2)$ are constants, satisfies (\ref
{EqPickandscondition}). See Huang, Zhang and Robeson \cite{HuangZR11} and Gneiting
\cite{Gneiting12} for further comments
on (\ref
{Defsin}) and more examples.

For $x=(x_1, \ldots, x_{N+1})\in\mathbb{S}^N$, its corresponding
spherical coordinate
$\theta=(\theta_1, \ldots, \theta_N)$ is defined as follows.
%
%
\begin{eqnarray}
\label{Eqsphericalcoordinate} x_1 &=& \cos\theta_1,\nonumber
\\
x_2 &=& \sin\theta_1 \cos\theta_2,\nonumber
\\
x_3 &=& \sin\theta_1 \sin\theta_2 \cos\theta_3,
\nonumber\\[-8pt]\\[-8pt]\nonumber
& \vdots&
\\
x_{N} &=& \sin\theta_1 \sin\theta_2 \cdots\sin
\theta_{N-1} \cos\theta_{N},\nonumber
\\
x_{N+1} &=& \sin\theta_1 \sin\theta_2 \cdots\sin
\theta_{N-1} \sin\theta_{N},\nonumber
\end{eqnarray}
where $0\leq\theta_i \leq\uppi$ for $1\leq i \leq N-1$ and $0 \leq
\theta_N < 2\uppi$.


We define the Gaussian field $\widetilde{X} = \{ \widetilde{X}(\theta
)\dvt \theta\in[0, \uppi]^{N-1}\times[0, 2 \uppi)\}$ by
$\widetilde{X}(\theta):= X(x)$ and denote by $\widetilde{C}$ the
covariance function of $\widetilde{X}$ accordingly.
The following elementary lemma characterizes the local behavior of the
spherical distance.
It provides a useful tool for establishing the relation between local
behaviors of covariance functions
$C$ and $\widetilde{C}$. Since we cannot find such a result in the
literature, for readers' convenience, we
provide here a short proof.

\begin{lemma}\label{Lemsphericaldomain}
Let $x, y\in\mathbb{S}^N$ and let $x$ be fixed. Then as $d(y,x) \to0$,
%
%
\begin{eqnarray}
\label{Eqequiv0} d^2(y,x) \sim(\varphi_1-
\theta_1)^2 + \bigl(\sin^2
\theta_1\bigr) (\varphi_2-\theta_2)^2+
\cdots+ \Biggl(\prod_{i=1}^{N-1}
\sin^2 \theta_i \Biggr) (\varphi_N-
\theta_N)^2.
\end{eqnarray}
Here and in the sequel, $\theta=(\theta_1, \ldots, \theta_N)$ and
$\varphi=(\varphi_1, \ldots, \varphi_N)$ are the spherical
coordinates of $x$ and $y$, respectively.
\end{lemma}

\begin{pf}
For $x, y\in\mathbb{S}^N$, we see that $d(y,x) \sim\llVert
y-x\rrVert $ as $d(y,x) \to0$, and
\begin{eqnarray*}
\llVert y-x\rrVert^2 &=& 2-2\cos(\varphi_1-
\theta_1) + 2(\sin\varphi_1\sin\theta_1)
\bigl[1-\cos(\varphi_2-\theta_2)\bigr]
\\
&&{} +\cdots+ 2 \Biggl(\prod_{i=1}^{N-1} \sin
\varphi_i\sin\theta_i \Biggr)\bigl[1-\cos(
\varphi_N-\theta_N)\bigr].
\end{eqnarray*}
It follows from the spherical coordinates that $d(y,x) \to0$ is
equivalent to $\llVert \varphi-\theta\rrVert \to0$. (There is an
exception for
$\theta$ with $\theta_N=0$, since for those $\varphi$ such that
$d(y,x) \to0$ and $\varphi_N$ tending to $2\uppi$, $\llVert \varphi
-\theta
\rrVert $ does not tend to $0$. In such case, we may treat $\theta_N$ as
$2\uppi$ instead of $0$ and this does not affect the result thanks to
the periodicity.) Therefore, as $d(y,x) \to0$, (\ref{Eqequiv0})
follows from Taylor's expansion.
\end{pf}

Next,\vspace*{1pt} we recall from Chan and Lai \cite{ChanL06} some results on the excursion
probability of Gaussian fields over
the Euclidean space.
Let $0<\alpha\leq2$ and let $\{W_t(s)\dvt s\in[0,\infty)^N\}$ ($t \in
\R^N$) be a family of Gaussian fields such that
%
%
\begin{eqnarray}
\label{EqlimitingGaussianfieldW} 
\E\bigl(W_t(s) \bigr)&=&-\llVert s\rrVert
^\alpha r_t\bigl(s/\llVert s\rrVert\bigr),\nonumber
\\
\operatorname{Cov} \bigl(W_t(s),W_t(v) \bigr)&=&\llVert s\rrVert
^\alpha r_t\bigl(s/\llVert s\rrVert\bigr)+\llVert v\rrVert
^\alpha r_t\bigl(v/\llVert v\rrVert\bigr)
\\
&&{} -\llVert s-v\rrVert^\alpha r_t\bigl((s-v)/\llVert s-v
\rrVert\bigr),\nonumber
\end{eqnarray}
where $r_t(\cdot)\dvtx \mathbb{S}^{N-1}\rightarrow\R_+$ is a continuous
function which satisfies
%
%
\begin{equation}
\label{Eqfunctionr} \sup_{v\in\mathbb{S}^{N-1} } \bigl\llvert r_t(v)-r_s(v)
\bigr\rrvert\to0\qquad\mbox{as } s\to t.
\end{equation}
Define
%
%
\begin{equation}
\label{EqgeneralPickandsconstant} H_\alpha^r(t)=\lim_{K\to\infty}K^{-N}
\int_0^\infty \mathrm{e}^u\P\Bigl\{ \sup
_{s\in[0,K]^N} W_t(s)\ge u \Bigr\}\,\mathrm{d}u.
\end{equation}
Denote by $H_\alpha$ the usual Pickands' constant, that is
\[
H_\alpha=\lim_{K\to\infty}K^{-N}\int
_0^\infty \mathrm{e}^u\P\Bigl\{\sup
_{s\in[0,K]^N} Z(s)\ge u \Bigr\}\,\mathrm{d}u,
\]
where $\{Z(s)\dvt s\in[0,\infty)^N\}$ is a Gaussian field such that
\begin{eqnarray*}
\E\bigl(Z(s) \bigr)=-\llVert s\rrVert^\alpha,\qquad \operatorname{Cov}
\bigl(Z(s),Z(v) \bigr)=\llVert s\rrVert^\alpha+ \llVert v\rrVert
^\alpha-\llVert s-v\rrVert^\alpha.
\end{eqnarray*}
It is clear that 
$H_\alpha^r(t)$ becomes $H_\alpha$ when $r_t\equiv1$.

Let $D\subset\R^N$ be a bounded $N$-dimensional Jordan measurable
set, that is, the boundary of $D$
has $N$-dimensional Lebesgue measure 0. Let $Y = \{Y(t), t \in\mathbb
R^N\}$
be a real-valued, centered Gaussian field such that its covariance
function $C_Y$ satisfies
%
%
\begin{equation}
\label{Eqlocallyisotropicproperty} C_Y(t,t+s)= 1-\llVert s\rrVert
^\alpha
r_t\bigl(s/\llVert s\rrVert\bigr) \bigl(1+\mathrm{o}(1)\bigr)\qquad\mbox{as
} \llVert s\rrVert\to0,
\end{equation}
for some constant $\alpha\in(0,2]$, uniformly over $t\in\bar{D}$,
the closure of $D$.

We will make use of the following theorem of Chan and Lai \cite{ChanL06}. One
can also apply similar results in Piterbarg
\cite{Piterbarg96}, Mikhaleva and Piterbarg \cite{MP96}, which are formulated
under somewhat different local
stationarity conditions. Having the functions $r_t(\cdot)$ in (\ref
{Eqlocallyisotropicproperty}) makes the
following theorem slightly easier to apply.

%
\begin{theorem}[(Chan and Lai \cite{ChanL06}, Theorem~2.1)]\label{ThmPickandsapproximation}
Let $D\subset\R^N$ be a bounded $N$-dimensional Jordan measurable
set. Suppose the Gaussian
field $\{Y(t)\dvt t\in\R^N\}$ satisfies condition (\ref
{Eqlocallyisotropicproperty}), in which
$r_t(\cdot)\dvtx \mathbb{S}^{N-1}\rightarrow\R_+$ is a continuous
function such that the convergence
(\ref{Eqfunctionr}) is uniform in $\bar{D}$ and $\sup_{t\in\bar
{D}, v\in\mathbb{S}^{N-1}} r_t(v)<\infty$. Then as $u\to\infty$,
\[
\P\Bigl\{\sup_{t\in D} Y(t) \ge u \Bigr\} \sim u^{2N/\alpha}
\Psi(u) \int_D H^r_\alpha(t)\,\mathrm{d}t.
\]
Here and in the sequel, $ \Psi(u) = (\sqrt{2 \uppi} u)^{-1} \mathrm{e}^{-u^2/2}$.
\end{theorem}

The lemma below establishes the relation between $H_\alpha^r(t)$ and
$H_\alpha$ for a special class of
functions $r_t(\cdot)$.

\begin{lemma}\label{LemgeneralandusualPickandsconstants}
Let $\{W_t(s)\dvt s\in[0,\infty)^N\}$ ($t \in\R^N$) be a family of
Gaussian fields satisfying (\ref{EqlimitingGaussianfieldW}) with
$r_t(v)=\llVert M_tv\rrVert ^\alpha$ for all $ v\in\mathbb{S}^{N-1}$,
where, for every $t \in\R^N$, $M_t$ is a non-degenerate $N\times N$
matrix. Then
$ H_\alpha^r(t) = \llvert \operatorname{det}M_t\rrvert H_\alpha
$ for each $t\in\R^N$.
\end{lemma}

\begin{pf} Let $t\in\R^N$ be fixed and consider the centered
Gaussian field $\overline{W}_t
= \{ \overline{W}_t(s), s\in[0,\infty)^N\}$ defined by $\overline
{W}_t(s) =W_t(M_t^{-1} s)$. Then by
(\ref{EqlimitingGaussianfieldW}), $\overline{W}_t$ satisfies
%
%
\begin{eqnarray}
\label{EqW2} \E\bigl( \overline{W}_t(s) \bigr) =-\llVert s\rrVert
^\alpha,\qquad \operatorname{Cov} \bigl(\overline{W}_t( s),
\overline{W}_t( v) \bigr)=\llVert s\rrVert^\alpha+\llVert v
\rrVert^\alpha-\llVert s-v\rrVert^\alpha.
\end{eqnarray}
Let $B_K=[0,K]^N$ and $M_tB_K=\{s\in\R^N\dvt \exists v\in B_K$ such that
$s=M_tv\}$.
Then $\operatorname{Vol}(M_tB_K)=\llvert \operatorname{det}M_t\rrvert
\operatorname{Vol}(B_K)$ and $\sup_{s\in
B_K} W_t(s) = \sup_{s\in M_tB_K}
\overline{W}_t(s)$, it follows from (\ref{EqgeneralPickandsconstant}) that
%
%
\begin{eqnarray}
\label{EqH-trans} H_\alpha^r(t)&=&\lim_{K\to\infty}
\frac{1}{\operatorname{Vol}(B_K)}\int_0^\infty \mathrm{e}^u\P
\Bigl\{\sup_{s\in B_K} W_t(s)\ge u \Bigr\}\,\mathrm{d}u\nonumber
\\
&=&\lim_{K\to\infty}\frac{\operatorname{Vol}(M_tB_K)}{\operatorname
{Vol}(B_K)}\frac
{1}{\operatorname{Vol}(M_tB_K)}\int
_0^\infty \mathrm{e}^u\P\Bigl\{\sup
_{s\in M_tB_K} \overline{W}_t(s)\ge u \Bigr\}\,\mathrm{d}u
\\
&=&\llvert\operatorname{det}M_t\rrvert\lim_{K\to\infty}
\frac{1}{\operatorname{Vol}(M_tB_K)}\int_0^\infty \mathrm{e}^u\P
\Bigl\{\sup_{s\in M_tB_K} \overline{W}_t(s)\ge u \Bigr\}\,\mathrm{d}u.\nonumber
\end{eqnarray}
Because of (\ref{EqW2}), we can modify the proofs in \cite
{QuallsW96} to show that
%
%
\begin{equation}
\label{EqH-orig} H_\alpha=\lim_{K\to\infty}\frac{1}{\operatorname
{Vol}(M_tB_K)}\int
_0^\infty \mathrm{e}^u\P\Bigl\{\sup
_{s\in M_tB_K} \overline{W}_t(s)\ge u \Bigr\}\,\mathrm{d}u.
\end{equation}
Comparing (\ref{EqH-trans}) and (\ref{EqH-orig}) gives the result.
\end{pf}

For any $T\subset\mathbb{S}^N$, we denote by $D \subset[0, \uppi
]^{N-1}\times[0, 2 \uppi)$ the set corresponding to $T$ under
the spherical coordinates (\ref{Eqsphericalcoordinate}).
We say that $T$ is an $N$-dimensional Jordan measurable set on $\mathbb
{S}^N$ if $D$ is an $N$-dimensional Jordan
measurable set in $\R^N$. Now we can prove our main result of this section.

\begin{theorem}\label{ThmPickandsapproximationsphere}
Let $\{X(x)\dvt x\in\mathbb{S}^N\}$ be a centered Gaussian random field
satisfying condition (\ref{EqPickandscondition})
and let $T\subset\mathbb{S}^N$ be an $N$-dimensional Jordan
measurable set on $\mathbb{S}^N$. Then as $u\to\infty$,
\[
\P\Bigl\{\sup_{x\in T} X(x) \ge u \Bigr\} \sim c^{N/\alpha
}\operatorname{Area}(T)H_\alpha u^{2N/\alpha} \Psi(u),
\]
where $\operatorname{Area}(T)$ denotes the spherical area of $T$ and
$c>0$ is
the constant in (\ref{EqPickandscondition}).
\end{theorem}

\begin{pf} For any $\theta\in[0, \uppi]^{N-1}\times[0, 2 \uppi)$, let
$M_\theta=c^{1/\alpha}\operatorname{diag}(1, \sin\theta_1, \ldots, \prod
_{i=1}^{N-1} \sin\theta_i)$
be the $N \times N$ diagonal matrix. 
By Lemma \ref{Lemsphericaldomain}, condition (\ref
{EqPickandscondition}) implies
\[
\widetilde{C}(\theta,\theta+\xi)= 1-\llVert\xi\rrVert^\alpha
r_\theta\bigl(\xi/\llVert\xi\rrVert\bigr) \bigl(1+\mathrm{o}(1)\bigr)\qquad
\mbox{as } \llVert\xi\rrVert\to0,
\]
where $r_\theta(\tau)=\llVert M_\theta\tau\rrVert ^\alpha$, $\forall
\tau\in
\mathbb{S}^{N-1}$. Then by Theorem \ref{ThmPickandsapproximation},
as $u\to\infty$,
%
%
\begin{equation}
\label{EqexcurionproboverTandD} \P\Bigl\{\sup_{x\in T} X(x) \ge u \Bigr
\}= \P
\Bigl\{\sup_{\theta
\in D} \widetilde{X}(\theta) \ge u \Bigr\} \sim
u^{2N/\alpha} \Psi(u) \int_D H^r_\alpha(
\theta)\,\mathrm{d}\theta.
\end{equation}
It follows from Lemma \ref{LemgeneralandusualPickandsconstants}
that for any $\theta\in[0, \uppi]^{N-1}\times[0, 2 \uppi)$
such that $M_\theta$ is non-degenerate (i.e., $\prod_{i=1}^{N-1} \sin
\theta_i\neq0$),
\[
H^r_\alpha(\theta)=c^{N/\alpha} \Biggl( \prod
_{i=1}^{N-1}\sin^{N-i} \theta_i
\Biggr)H_\alpha. %
\]
Note that $( \prod_{i=1}^{N-1}\sin^{N-i} \theta_i )\,\mathrm{d}\theta$ is the
spherical area element and $M_\theta$ is non-degenerate for almost every
$\theta\in D$, we obtain
\[
\int_D H^r_\alpha(\theta)\,\mathrm{d}\theta=
c^{N/\alpha}\operatorname{Area}(T)H_\alpha.
\]
Plugging this into (\ref{EqexcurionproboverTandD}) gives the
desired result.
\end{pf}


\subsection{Standardized spherical fractional Brownian motion}\label{secSSFBM}
Theorem \ref{ThmPickandsapproximationsphere} 
provides a nice approximation to the excursion probability for locally
isotropic Gaussian random fields on $\mathbb S^N$ whose covariance functions
satisfy (\ref{EqPickandscondition}). When the local behavior of the
covariance function becomes more complicated, Theorem \ref
{ThmPickandsapproximationsphere} may not be applicable anymore.
However, we can still apply Lemma \ref{Lemsphericaldomain} to find
the corresponding local behavior of covariance function under spherical
coordinates and then apply Theorem \ref{ThmPickandsapproximation} to
obtain the asymptotics for the excursion probability. In the following,
we use spherical fractional Brownian motion on sphere as an
illustrating example.

Let $o$ be a fixed point on $\mathbb{S}^N$. The spherical fractional
Brownian motion (SFBM) $B_\beta= \{B_\beta(x)\dvt x \in\mathbb S^N\}$
is defined by Istas \cite{Istas05} as a centered real-valued Gaussian random
field such that $B_\beta(o)=0 $ and
\[
\E\bigl(B_\beta(x)-B_\beta(y)\bigr)^2 =
d^{2\beta}(x,y) \qquad\forall x,y\in\mathbb{S}^N,
\]
where $\beta\in(0, 1/2]$. It follows immediately that
\[
\operatorname{Cov}\bigl(B_\beta(x), B_\beta(y)\bigr)=\tfrac{1}{2}
\bigl(d^{2\beta
}(x,o)+d^{2\beta}(y,o)-d^{2\beta}(x,y) \bigr).
\]
Without loss of generality, we take $o=(1,0,\ldots,0)\in\R^{N+1}$,
whose corresponding spherical coordinate
is $(0,\ldots,0)\in\R^N$. We consider the standardized SFBM $X = \{
X(x)\dvt x\in\mathbb{S}^N \setminus\{o\}\}$
defined by
%
%
\begin{equation}
\label{EqX3} X(x)=\frac{B_\beta(x)}{d^\beta(x,o)} \qquad\forall x\in
\mathbb
{S}^N \setminus\{o\}.
\end{equation}
Then the covariance of $X$ is
\[
C(x,y)=\operatorname{Cov}\bigl(X(x), X(y)\bigr)=\frac{d^{2\beta
}(x,o)+d^{2\beta
}(y,o)-d^{2\beta}(x,y)}{2d^\beta(x,o)d^\beta(y,o)}.
\]
Note that, under the spherical coordinates, $d(x,o)=\theta_1$ and
$d(y,o)=\varphi_1$, together with
Lemma~\ref{Lemsphericaldomain}, we obtain that the covariance
function of the corresponding Gaussian field
$\widetilde{X}$ satisfies
\begin{eqnarray*}
\widetilde{C}(\theta,\varphi)&=&\operatorname{Cov}\bigl(\widetilde
{X}(\theta),
\widetilde{X}(\varphi)\bigr)
\\
&=&1-\bigl(1+\mathrm{o}(1)\bigr)
\\
&&\hspace*{15pt}{}\times \frac{1}{2\theta
_1^{2\beta}} \Biggl[(
\varphi_1-\theta_1)^2 + \bigl(
\sin^2 \theta_1\bigr) (\varphi_2-\theta
_2)^2
 +\cdots+ \Biggl(\prod_{i=1}^{N-1}
\sin^2 \theta_i \Biggr) (\varphi_N-
\theta_N)^2 \Biggr]^\beta
\end{eqnarray*}
as $d(x,y)\to0$. Let
\begin{eqnarray*}
M_\theta&=&\frac{1}{2^{1/(2\beta)}\theta_1}\operatorname{diag} \Biggl(1,
\sin
\theta_1, \ldots, \prod_{i=1}^{N-1}
\sin\theta_i \Biggr),
\\
r_\theta(\tau)&=&\llVert M_\theta\tau\rrVert^{2\beta}
\qquad\forall\tau\in\mathbb{S}^{N-1},
\end{eqnarray*}
and $\xi=\varphi-\theta$, then as $\llVert \xi\rrVert \to0$,
\[
\widetilde{C}(\theta,\theta+\xi)=1-\llVert\xi\rrVert^{2\beta}
r_\theta\bigl(\xi/\llVert\xi\rrVert\bigr) \bigl(1+\mathrm{o}(1)\bigr).
\]
Let $T\subset\mathbb{S}^N$ be an $N$-dimensional Jordan measurable
set such that $o\notin\bar{T}$,
and denote its corresponding domain under the spherical coordinates by
$D$, which implies $\theta_1 \neq0$
for any $\theta\in\bar{D}$. By Theorem \ref{ThmPickandsapproximation},
as $u\to\infty$,
\[
\P\Bigl\{\sup_{x\in T} X(x) \ge u \Bigr\}= \P\Bigl\{\sup
_{\theta
\in D} \widetilde{X}(\theta) \ge u \Bigr\} \sim u^{N/\beta}
\Psi(u) \int_D H^r_{2\beta}(\theta)\,\mathrm{d}
\theta.
\]
For any $\theta$ such
that $M_\theta$ is non-degenerate (i.e., $\prod_{i=1}^{N-1} \sin
\theta_i\neq0$), Lemma \ref{LemgeneralandusualPickandsconstants}
gives
\[
H^r_{2\beta}(\theta)=\frac{1}{2^{N/(2\beta)}\theta_1^{N}} \Biggl( \prod
_{i=1}^{N-1}\sin^{N-i} \theta_i
\Biggr)H_{2\beta}. %
\]
Therefore, as $u\to\infty$,
%
%
\begin{eqnarray}
\label{EqX4} \P\Bigl\{\sup_{x\in T} X(x) \ge u \Bigr\} &\sim&
u^{N/\beta} \Psi(u) 2^{-N/(2\beta)} H_{2\beta} \int
_D \theta_1^{-N} \Biggl( \prod
_{i=1}^{N-1}\sin^{N-i}
\theta_i \Biggr) \,\mathrm{d}\theta.
\end{eqnarray}

%
\begin{remark}\label{Re;25}
Comparing the excursion probabilities in (\ref{EqX4}) for the
standardized SFBM $X$ and in Theorem
\ref{ThmPickandsapproximationsphere} for the isotropic Gaussian
field $W^\beta$, which is defined in Estrade and Istas \cite{EsIstas10}, we see
that the constant in (\ref{EqX4}) is more complicated.
\end{remark}

\section{Smooth isotropic Gaussian fields on sphere}
In this section, we study the excursion probability of smooth isotropic
Gaussian fields on sphere.
Related to the results in this section, we mention that \cite{ChengS13}
have determined the height distribution and overshoot distribution of
local maxima of smooth isotropic Gaussian
random fields on sphere.

\subsection{Preliminaries}
Given $\la>0$ and an integer $n \ge0$, the \emph{ultraspherical
polynomial} (or \emph{Gegenbauer polynomial}) of degree $n$,
denoted by $P_n^\la(t)$, is defined by
the expansion
\[
\bigl(1-2rt +r^2\bigr)^{-\la}=\sum
_{n=0}^\infty r^n P_n^\la(t),
\qquad t\in[-1,1].
\]
For $\la=0$, we follow Schoenberg \cite{Schoenberg42} and define
$P_n^0(t)=\cos
(n\arccos t)=T_n(t)$, where $T_n$ ($n\ge0$)
are \emph{the Chebyshev polynomials of the first kind} defined by the expansion
\[
\frac{1-rt}{1-2rt +r^2}=\sum_{n=0}^\infty
r^n T_n(t), \qquad t\in[-1,1].
\]

For reference later on, we recall the following formulae on $P_n^\la$.

\begin{longlist}[(ii)]
\item[(i)] For all $n\ge0$, $P_n^0(1)=1$, and if $\la>0$ (cf. Szeg{\H{o}} \cite{Szego75}, page~80),
%
%
\begin{equation}
\label{Equspolynomialsat1} P_n^\la(1)=\pmatrix{n+2\la-1
\cr
n}.
\end{equation}

\item[(ii)] For all $n\ge0$,
%
%
\begin{equation}
\label{Eqdifferentiateuspolynomial0} \frac{\mathrm{d}}{\mathrm{d}t}P_n^0(t) = n
P_{n-1}^1(t),
\end{equation}
and if $\la>0$ (cf. Szeg{\H{o}} \cite{Szego75}, page~81),
%
%
\begin{equation}
\label{Eqdifferentiateuspolynomial} \frac{\mathrm{d}}{\mathrm{d}t}P_n^\la(t) = 2\la
P_{n-1}^{\la+1}(t).
\end{equation}

The following theorem by Schoenberg \cite{Schoenberg42} characterizes the covariance
function of an isotropic
Gaussian field on sphere (see also Gneiting \cite{Gneiting12}).
\end{longlist}

%
\begin{theorem}\label{ThmSchoenberg} Let $N\ge1$, then a continuous
function $C(\cdot, \cdot)\dvtx
\mathbb{S}^N\times\mathbb{S}^N \rightarrow\R$ is the covariance of
an isotropic Gaussian field
on $\mathbb{S}^N$ if and only if it has the form
\[
C(x,y)= \sum_{n=0}^\infty a_n
P_n^\la\bigl({\langle}x, y \rangle\bigr), \qquad x, y \in
\mathbb{S}^N,
\]
where $\la=(N-1)/2$, $a_n \geq0$ and $\sum_{n=0}^\infty a_nP_n^\la
(1) <\infty$.
\end{theorem}

%
\begin{remark}
Note that for the case of $N=1$ and $\la=0$, $\sum_{n=0}^\infty
a_nP_n^0(1) <\infty$ is equivalent to $\sum_{n=0}^\infty a_n <\infty$;
while for $N\ge2$ and $\la=(N-1)/2$, (\ref{Equspolynomialsat1})
implies that $\sum_{n=0}^\infty a_nP_n^\la(1) <\infty$
is equivalent to $\sum_{n=0}^\infty n^{N-2}a_n <\infty$.
\end{remark}

When $N=2$ and $\la=1/2$, $P_n^\la$ ($n \ge0$) become the \emph
{Legendre polynomials}. For more
results on isotropic Gaussian fields on $\mathbb{S}^2$, we refer to
Marinucci and Peccati \cite{MarinucciPecati07}. Regularity and smoothness properties of Gaussian
field $\{X(x)\dvt x \in\mathbb{S}^2\}$
have recently been obtained by Lang and Schwab \cite{LS13} in terms of the
corresponding angular power spectrum.

The following statement $({\mathbf A}1)$ is a smoothness condition for
Gaussian fields on sphere. In Lemma \ref{LemC^3sphere} below, we
show that it implies $X(\cdot)\in C^2(\mathbb{S}^N)$ a.s.
\begin{longlist}
\item[$({\mathbf A}1)$.] The covariance $C(\cdot, \cdot)$ of $\{
X(x)\dvt
x\in\mathbb{S}^N\}$ satisfies
\[
C(x,y)= \sum_{n=0}^\infty a_n
P_n^\la\bigl({\langle}x, y \rangle\bigr), \qquad x, y \in
\mathbb{S}^N,
\]
where $\la= \frac{N-1} 2$, $a_n \geq0$ and $\sum_{n=1}^\infty
n^{N+8}a_n<\infty$ if $N\ge2$; $\sum_{n=1}^\infty n^{10}a_n<\infty$
if $N=1$.
\end{longlist}

%
\begin{lemma}\label{LemC^3sphere}
Let $\{X(x)\dvt x\in\mathbb{S}^N\}$ be an isotropic Gaussian field such
that $({\mathbf A}1)$ is fulfilled. Then $X(\cdot)\in C^2(\mathbb
{S}^N)$ a.s.
\end{lemma}

\begin{pf} We first consider $N\ge2$. By Theorem \ref
{ThmSchoenberg}, each $P_n^\la({\langle}t, s \rangle)$ is the
covariance of an isotropic Gaussian field on $\mathbb{S}^N$ and hence
the Cauchy--Schwarz inequality implies
%
%
\begin{equation}
\label{Eqbounduspolynomials} \bigl\llvert P_n^\la\bigl({\langle}x, y
\rangle\bigr)\bigr\rrvert\le P_n^\la\bigl({\langle}x, x
\rangle\bigr) = P_n^\la(1) \qquad\forall x,y\in
\mathbb{S}^N.
\end{equation}
Combining $({\mathbf A}1)$ with (\ref{Equspolynomialsat1}), (\ref
{Eqdifferentiateuspolynomial}) and (\ref{Eqbounduspolynomials}),
together with the fact $P_{0}^{\la}(t)\equiv1$, we obtain that there
exist positive constants $M_1$ and $M_2$ such that
\[
\sup_{t\in[-1,1]}\sum_{n=0}^\infty
a_n \biggl\llvert\biggl(\frac
{\mathrm{d}^5}{\mathrm{d}t^5}P_n^\la(t)
\biggr)\biggr\rrvert\le M_1\sum_{n=5}^\infty
a_n P_{n-5}^{\la+5}(1)\le M_2 \sum
_{n=1}^\infty n^{N+8}a_n<
\infty. %
\]
This shows that $C(\cdot, \cdot) \in C^5(\mathbb{S}^N\times\mathbb
{S}^N)$. The proof for $N=1$ is similar once we apply both (\ref
{Eqdifferentiateuspolynomial0}) and (\ref
{Eqdifferentiateuspolynomial}). Therefore, by arguments via charts (cf.
Auffinger \cite{AuffingerPhD})
and the results in Potthoff \cite{Potthoff10} (though the results therein
are for
$X(\cdot)\in C^1$, they can be extended easily to the case of
higher-order smoothness), we conclude that $X(\cdot)\in C^2(\mathbb
{S}^N)$ a.s.
\end{pf}

By Schoenberg \cite{Schoenberg42} or Gneiting \cite{Gneiting12}, $C(\cdot, \cdot
)$ is a
covariance function on $\mathbb{S}^N$ for every $N\ge1$
if and only if it has the form
\[
C(x,y)= \sum_{n=0}^\infty b_n {
\langle}x, y \rangle^n, \qquad x, y \in\mathbb{S}^N,
\]
where $b_n\ge0$ and $\sum_{n=0}^\infty b_n<\infty$. Then similarly
to $({\mathbf A}1)$, we may state the smoothness condition~$({\mathbf A}1')$ below for this special class of Gaussian fields on sphere.
\begin{longlist}
\item[$({\mathbf A}1')$.] The covariance $C(\cdot, \cdot)$ of $\{
X(x)\dvt
x\in\mathbb{S}^N\}$ satisfies
\[
C(x,y)= \sum_{n=0}^\infty b_n {
\langle}x, y \rangle^n, \qquad x, y \in\mathbb{S}^N,
\]
where $b_n\ge0$ and $\sum_{n=0}^\infty n^5 b_n<\infty$.
\end{longlist}

We obtain below an analogue of Lemma \ref{LemC^3sphere}. Since the
proof is similar, it is omitted.

\begin{lemma}\label{LemC^3sphereinfinity}
Let $\{X(x)\dvt x\in\mathbb{S}^N\}$ be an isotropic Gaussian field such
that $({\mathbf A}1')$ is fulfilled. Then $X(\cdot)\in C^2(\mathbb
{S}^N)$ a.s.
\end{lemma}

\subsection{Excursion probability}
Let $\chi(A_u(X,\mathbb{S}^N))$ be the Euler characteristic of
excursion set $A_u(X,\mathbb{S}^N)
= \{x\in\mathbb{S}^N\dvt X(x)\geq u\}$ (cf. Adler and Taylor \cite{AT07}).
Denote by $H_j(x)$ the Hermite
polynomial of order $j$, that is,
\[
H_j(x) = (-1)^j \mathrm{e}^{x^2/2} \frac{\mathrm{d}^j}{\mathrm{d}x^j}
\bigl( \mathrm{e}^{-x^2/2} \bigr).
\]
Denote $\omega_j = \frac{2\uppi^{(j+1)/2}}{\Gamma((j+1)/{2})}$, the
spherical area of the $j$-dimensional unit sphere $\mathbb{S}^j$.

Before stating our results, we need another regularity condition for
the Gaussian field.
\begin{longlist}
\item[$({\mathbf A}2)$.] For each $x\in\mathbb{S}^N$, the joint
distribution of
$(X(x), \nabla X(x), \nabla^2X(x))$ is non-degenerate.
\end{longlist}

%
\begin{lemma}\label{LemMECsphere} Let $\{X(x)\dvt x\in\mathbb{S}^N\}$
be a centered, unit-variance, isotropic Gaussian
field satisfying $({\mathbf A}1)$ and $({\mathbf A}2)$. Then
\begin{eqnarray*}
\E\bigl\{\chi\bigl(A_u\bigl(X,\mathbb{S}^N\bigr)\bigr)
\bigr\} = \sum_{j=0}^N
\bigl(C'\bigr)^{j/2} \mathcal{L}_j \bigl(
\mathbb{S}^N\bigr) \rho_j(u),
\end{eqnarray*}
where the constant $C'$ is defined as
%
%
\begin{eqnarray}
\label{DefC} C'= \cases{ \displaystyle (N-1)\sum_{n=1}^\infty
\pmatrix{n+N-1
\cr
N}a_n, &\quad if $N\ge2$,
\cr
\displaystyle\sum
_{n=1}^\infty n^2 a_n, &\quad if
$N=1$,}
\end{eqnarray}
and where $\rho_0(u)= (2\uppi)^{-1/2}\int_u^\infty \mathrm{e}^{-x^2/2}\,\mathrm{d}x$,
$\rho_j(u) = (2\uppi)^{-(j+1)/2} H_{j-1} (u) \mathrm{e}^{-u^2/2}$
for $j\geq1$ and
%
%
\begin{eqnarray}
\label{EqL-Kcurvature} \mathcal{L}_j \bigl(\mathbb{S}^N\bigr) =
\cases{\displaystyle 2 \pmatrix{N
\cr
j}\displaystyle\frac{\omega_N}{\omega_{N-j}}, &\quad if $N-j$ is even,
\cr
0, &
\quad otherwise}
\end{eqnarray}
(for $j=0, 1, \ldots, N$) are the Lipschitz--Killing curvatures of
$\mathbb{S}^N$ (cf. (6.3.8) in Adler and Taylor~\cite{AT07}).
\end{lemma}

%
\begin{remark}
In Lemma \ref{LemMECsphere}, if condition $({\mathbf A}1)$ is replaced
by $({\mathbf A}1')$, then it can be seen from
the proof below that the result still holds with $C'$ being replaced by
$C'=\sum_{n=1}^\infty n b_n$.
\end{remark}

\begin{pf*}{Proof of Lemma \ref{LemMECsphere}} By Theorem 12.4.1 in Adler and Taylor
\cite{AT07},
we only need to show that the Lipschitz--Killing curvatures induced by
$X$ on $\mathbb{S}^N$ are
$\mathcal{L}_j (X,\mathbb{S}^N)=(C')^{j/2}\mathcal{L}_j (\mathbb
{S}^N)$ for $j = 0, 1, \ldots, N$.

The Riemannian structure induced by $X$ on $\mathbb{S}^N$ is defined as
\[
g_{x_0}^{X, \mathbb{S}^N}(\xi_{x_0}, \sigma_{x_0}):= \E
\bigl\{(\xi_{x_0} X)\cdot(\sigma_{x_0} X)\bigr\}=
\xi_{x_0}\sigma_{x_0} C(x,y)\mid_{x=y=x_0} \qquad\forall
x_0 \in\mathbb{S}^N,
\]
where\vspace*{1pt} $\xi_{x_0}, \sigma_{x_0} \in T_{x_0}\mathbb{S}^N$, the tangent
space of $\mathbb{S}^N$ at $x_0$
(cf. Adler and Taylor \cite{AT07}, page~305). We may choose two
smooth curves
on $\mathbb{S}^N$, say $\gamma(t)$, $\tau(s)$, $t, s\in[0,1]$, such
that $\gamma(0)=\tau(0)=x_0$ and $\gamma'(0)=\xi_{x_0}, \tau
'(0)=\sigma_{x_0}$. We first consider $N\ge2$, then
\begin{eqnarray*}
\xi_{x_0}\sigma_{x_0} C(x,y)| _{x=y=x_0} &=&
\frac{\partial}{\partial
t}\frac{\partial}{\partial s} C\bigl(\gamma(t), \tau(s)\bigr)\Big|
_{t=s=0}
\\
&=&\frac{\partial}{\partial t}\frac{\partial}{\partial s} \sum
_{n=0}^\infty
a_n P_n^\la\bigl({\bigl\langle}\gamma(t),
\tau(s) \bigr\rangle\bigr)\Big|_{t=s=0}
\\
&=& \frac{\partial}{\partial t}\sum_{n=1}^\infty
a_n (N-1) P_{n-1}^{\la+1}\bigl({\bigl\langle}
\gamma(t), x_0 \bigr\rangle\bigr) {\bigl\langle}\gamma(t),
\sigma_{x_0} \bigr\rangle\Big| _{t=0}
\\
&=& \sum_{n=2}^\infty a_n (N-1)
(N+2) P_{n-2}^{\la+2}\bigl({\langle}x_0,
x_0 \rangle\bigr) {\langle}\xi_{x_0}, x_0
\rangle{\langle}x_0, \sigma_{x_0} \rangle
\\
&&{} + \sum_{n=1}^\infty a_n
(N-1) P_{n-1}^{\la+1}\bigl( {\langle}x_0,
x_0 \rangle\bigr) {\langle}\xi_{x_0}, \sigma_{x_0}
\rangle
\\
&=& \Biggl(\sum_{n=1}^\infty a_n
(N-1) P_{n-1}^{\la+1}(1) \Biggr) {\langle}\xi_{x_0},
\sigma_{x_0} \rangle= C'{\langle}\xi_{x_0},
\sigma_{x_0} \rangle,
\end{eqnarray*}
where the third and fourth equalities follow from (\ref
{Eqdifferentiateuspolynomial}), while the fifth equality is due to
the facts ${\langle}x_0, x_0 \rangle=1$ and ${\langle}\xi_{x_0},
x_0 \rangle= {\langle}\sigma_{x_0}, x_0 \rangle=0$, since the
vector $x_0$ is always orthogonal to its tangent space. The case $N=1$
can be proved similarly once we apply (\ref
{Eqdifferentiateuspolynomial0}) instead of (\ref{Eqdifferentiateuspolynomial}).

Hence the induced metric is
\[
g_{x_0}^{X, \mathbb{S}^N}(\xi_{x_0}, \sigma_{x_0}) =
C' {\langle}\xi_{x_0}, \sigma_{x_0} \rangle
\qquad\forall x_0 \in\mathbb{S}^N.
\]
By the definition of Lipschitz--Killing curvatures, one has
$
\mathcal{L}_j (X,\mathbb{S}^N) = (C')^{j/2} \mathcal{L}_j (\mathbb{S}^N)$,
where $\mathcal{L}_j (\mathbb{S}^N)$ are the original
Lipschitz--Killing curvatures of $\mathbb{S}^N$ given by (\ref{EqL-Kcurvature}).
We have finished the proof.
\end{pf*}

Applying Lemma \ref{LemMECsphere} and Theorem 14.3.3 in Adler and Taylor \cite
{AT07}, we obtain
immediately the following approximation for the excursion probability.

\begin{theorem}\label{ThmEulerapproximation}
Suppose the conditions in Lemma \ref{LemMECsphere} hold. Then, under
the notation therein,
there exists a constant $\alpha_0>0$ such that as $u\to\infty$,
%
%
\begin{equation}
\label{EqEulerapproximation} \P\Bigl\{\sup_{x\in\mathbb{S}^N} X(x) \ge
u \Bigr\} = \sum
_{j=0}^N \bigl(C'
\bigr)^{j/2} \mathcal{L}_j \bigl(\mathbb{S}^N
\bigr) \rho_j(u)+ \mathrm{o}\bigl(\mathrm{e}^{-\alpha_0
u^2 - u^2/2}\bigr).
\end{equation}
\end{theorem}

%
\begin{remark} The following are some remarks.
\begin{itemize}
\item
Under the conditions in Theorem \ref{ThmEulerapproximation}, the
covariance function $C$ satisfies (\ref{EqPickandscondition}) with
$\alpha=2$.
Since for $\alpha=2$, Pickands' constant $H_2=\uppi^{-N/2}$, one can
check that the approximation in Theorem \ref
{ThmPickandsapproximationsphere} only provides the leading term of the
approximation in Theorem \ref{ThmEulerapproximation}. This also
affects the errors in two approximations: the error in the former one
is only $\mathrm{o}(u^{N-1}\mathrm{e}^{-u^2/2})$, while the error in the latter one is
$\mathrm{o}(\mathrm{e}^{-\alpha_0 u^2-u^2/2})$.
\item By applying the tube method, Sun \cite{Sun93} gave a two-term
approximation formula for the excursion probability of a class of differentiable
Gaussian random field $\{X(x), x \in I\}$, where $I \subset\R^N$ is a
bounded convex set. Her results can be applied to provide a two-term
approximation for the excursion probability in (\ref
{EqEulerapproximation}) for some special cases. See
Park and Sun \cite{ParkSun98}, page~73.
\item Recently Marinucci and Vadlamani \cite{MarinucciV13} have computed the
Lipschitz--Killing
curvatures of excursion set and derived a very precise approximation
for the excursion probability of a class of nonlinear functionals of a
smooth Gaussian random field on $\mathbb S^2$. In the linear case
(i.e., $q=1$) Theorem 21 of Marinucci and Vadlamani \cite{MarinucciV13} is a special
case of (\ref{EqEulerapproximation}) with $N=2$.
\end{itemize}
\end{remark}

If the sphere $\mathbb{S}^N$ is replaced by a more general subset $T
\subset\mathbb{S}^N$, by revising Lemma \ref{LemMECsphere} and
applying Theorem 14.3.3 in Adler and Taylor \cite{AT07} again, we obtain the
following corollary.

\begin{corollary}\label{CorEulerapproximation}
Suppose the conditions in Lemma \ref{LemMECsphere} hold. Let
$T\subset\mathbb{S}^N$ be a $k$-dimensional, locally convex, regular
stratified manifold (cf. Adler and Taylor \cite{AT07}, page~198),
then there
exists $\alpha_0>0$ such that as $u\to\infty$,
%
%
\begin{equation}
\P\Bigl\{\sup_{x\in T} X(x) \ge u \Bigr\} = \sum
_{j=0}^k \bigl(C'
\bigr)^{j/2} \mathcal{L}_j (T) \rho_j(u)+ \mathrm{o}\bigl(\mathrm{e}^{-\alpha_0 u^2 - u^2/2}\bigr),
\end{equation}
where $\mathcal{L}_j (T)$ are the Lipschitz--Killing curvatures of $T$
(cf. Adler and Taylor \cite{AT07}, page~175), $C'$~and $\rho_j(u)$
are as in
Lemma \ref{LemMECsphere}.
\end{corollary}

The parameter set $T\subset\mathbb{S}^N$ in Corollary \ref
{CorEulerapproximation} is assumed to be nice enough.
Roughly speaking, it looks like a convex set and can be decomposed into
several smooth manifolds, see Adler and Taylor \cite{AT07} for a rigorous
definition. Also, the $j$th Lipschitz--Killing
curvature $\mathcal{L}_j (T)$ can be viewed as the measure of the
$j$-dimensional boundary of $T$. One may use Steiner's formula
(Adler and Taylor \cite{AT07}, page~142) to compute the
Lipschitz--Killing curvatures of $T$ exactly. In particular, if $T$ is a
semisphere of dimension one, then
$\mathcal{L}_0 (T)=1$ and $\mathcal{L}_1 (T)=\uppi$. If $T$ is a
semisphere of dimension two, then $\mathcal{L}_0 (T)
=1$, $\mathcal{L}_1 (T)=\uppi$ and $\mathcal{L}_2 (T)=2\uppi$. More
generally, if $T$ is a $k$-dimensional, locally
convex, regular stratified manifold, then $\mathcal{L}_0 (T)$ is the
Euler characteristic,
$\mathcal{L}_k (T)$ is the volume and $\mathcal{L}_{k-1} (T)$ is half
of the surface area. For the other
$\mathcal{L}_j (T)$, $1\le j\le k-2$, we can apply Steiner's formula
to find their values.

Lastly, to further illustrate the main results of this paper, we give
more examples on approximating
the excursion probability of Gaussian fields on spheres, including both
smooth and non-smooth cases.

\begin{example}
The canonical Gaussian field on $\mathbb{S}^N$, denoted by $X$, has
covariance function given by
$C(x,y)={\langle}x, y \rangle$ (cf. Adler and Taylor \cite{AT07}). Since
$C(x,y)=\cos d(x,y)$, it satisfies
\[
C(x,y)= 1-\tfrac{1}{2} d^2(x,y) \bigl(1+\mathrm{o}(1)\bigr),\qquad
\mbox{as } d (x,y)\to0.
\]
Applying Theorem \ref{ThmPickandsapproximationsphere} with
$T=\mathbb{S}^N$, $c=1/2$ and $\alpha=2$, we obtain an approximation
to the excursion probability:
\[
\P\Bigl\{\sup_{x\in\mathbb{S}^N} X(x) \ge u \Bigr\} \sim
2^{-N/2}\operatorname{Area}\bigl(\mathbb{S}^N\bigr)H_2 u^N
\Psi(u) = (2\uppi)^{-(N+1)/2}\omega_N u^{N-1}\mathrm{e}^{-u^2/2}.
\]
However, by applying Theorem \ref{ThmEulerapproximation} with
$C'=1$, we get a more precise approximation:
\[
\P\Bigl\{\sup_{x\in\mathbb{S}^N} X(x) \ge u \Bigr\} = \sum
_{j=0}^N \mathcal{L}_j \bigl(
\mathbb{S}^N\bigr) \rho_j(u)+ \mathrm{o}\bigl(\mathrm{e}^{-\alpha_0 u^2 - u^2/2}
\bigr).
\]
\end{example}

\begin{example} Consider the Hamiltonian of the pure $p$-spin model on
$\mathbb{S}^{N-1}$
\[
H_{N,p}(x)=\frac{1}{N^{(p-1)/2}}\sum_{i_1,\ldots, i_p=1}^N
J_{i_1,\ldots, i_p}x_{i_1}\cdots x_{i_p} \qquad\forall
x=(x_1,\ldots, x_N)\in\mathbb{S}^{N-1},
\]
where $J_{i_1,\ldots, i_p}$ are independent standard Gaussian random
variables. Then $H_{N,p}$ and $H_{N,p'}$
are independent for any $p\neq p'$ and
\[
\E\bigl\{H_{N,p}(x)H_{N,p}(y)\bigr\}=\frac{1}{N^{p-1}}{
\langle}x, y \rangle^p.
\]
Let $(b_p)_{p\geq2}$ be a sequence of positive numbers such that $\sum
_{p=2}^\infty2^p b_p<\infty$ and define
\[
X(x)=\sum_{p=2}^\infty b_p
H_{N,p}(x).
\]
Then $X$ is a smooth Gaussian random field on $\mathbb{S}^{N-1}$ with
covariance
\[
C(x,y)=\sum_{p=2}^\infty\frac{b_p^2}{N^{p-1}}{
\langle}x, y \rangle^p.
\]
We can apply Theorem \ref{ThmEulerapproximation} or Corollary \ref
{CorEulerapproximation} to approximate the excursion probability.
\end{example}

\begin{example} Consider the Gaussian field $\{X(x)\dvt x\in\mathbb
{S}^N\}$ with covariance structure $C(x,y)=1-\frac{2}{\uppi} d(x,y)$
(cf. Zuo \cite{Zuo03}, Remark 3.3). Since
$d(x,y)=\arccos{\langle}x,
y\rangle$, we have
%
%
\begin{eqnarray}
C(x,y)
= \sum_{n=0}^\infty
\frac{(2n)!}{4^n(n!)^2(2n+1)}{\langle}x, y \rangle^{2n+1}:=\sum
_{n=0}^\infty b_n {\langle}x, y \rangle
^n.
\end{eqnarray}
It is easy to check that $\sum_{n=0}^\infty n b_n =\infty$, $({\mathbf
A}1')$ is not satisfied and hence Theorem \ref{ThmEulerapproximation}
is not applicable. Instead, we may use Theorem \ref
{ThmPickandsapproximation} to get an approximation to the excursion probability.
This result allows one to construct confidence regions for the true
projection median defined in Zuo (\cite{Zuo03}, Section~3)
without using the bootstrapping techniques.
\end{example}


\section*{Acknowledgements}
The authors thank Professor Enkelejd Hashorva,
Dr. Lanpeng Ji for stimulating
discussions and Dr. Xiaohui Liu for pointing out the connection between
Gaussian random fields on sphere
and projection depth functions in Zuo \cite{Zuo03}. They thank the referees
for their constructive comments which
have led to improvements of the
manuscript.

Research partially supported by NSF Grants DMS-13-09856 and DMS-13-07470.


%

\printhistory
\end{document}